\documentclass[12pt,a4paper]{amsart}
\usepackage{latexsym}
\usepackage{amsmath}
\usepackage{amsthm}
\usepackage{amssymb}
\usepackage{amscd}
\usepackage[all]{xy}
\usepackage[cp850]{inputenc}
\usepackage[mathscr]{eucal}
\usepackage{color}
\theoremstyle{plain}
\newcommand{\adef}{\begin{defin}}
\newcommand{\zdef}{\end{defin}}
\newcommand{\To}{\longrightarrow}

\newtheorem{theorem}{Theorem}[section]

\newtheorem{defin}[theorem]{Definition}
\newtheorem{prop}[theorem]{Proposition}
\newtheorem{corollary}[theorem]{Corollary}
\newtheorem{lemma}[theorem]{Lemma}
\theoremstyle{remark}

\newcommand{\lop}{\curvearrowright }

\newcommand{\N}{\mathbb{N}}

\newcommand{\e}{\varepsilon}

\newcommand{\Dom}{\mathrm{Dom}\,}
\newcommand{\Ran}{\mathrm{Ran}\,}
\newcommand{\KP}{{\sf K}\hspace{-1pt}{\sf P}}

\newcommand{\ra}{\rightarrow}

\newcommand{\dist}{\textrm{dist}}
\newcommand{\im}{\textrm{Im}\,}

\title{Quasilinear duality and inversion in Banach spaces}

\author{Jes\'{u}s M. F. Castillo}
\address{Universidad de Extremadura,\\ Instituto de Matem\'aticas Imuex\\
Avenida de Elvas\\ 06071-Badajoz\\ Spain} \email{castillo@unex.es}

\author{Manuel Gonz\'alez}
\address{Departamento de Matem\'aticas, Universidad de Cantabria, Avenida de los Castros s/n, 39071 Santander, Spain}
\email{manuel.gonzalez@unican.es}

\thanks{This research was supported in part by MINCIN Project PID2019-103961GB-C21. The research of the first author was supported in part by project IB20038 de la Junta de Extremadura}

\begin{document}

\begin{abstract} We present a unified approach to the processes of inversion and duality for quasilinear and $1$-quasilinear maps; in particular, for centralizers and differentials generated by interpolation methods.\end{abstract}

\maketitle

\thispagestyle{empty}

\section{Introduction}
In this paper we present a unified approach to the inversion and duality phenomena for centralizers and differential maps generated by interpolation processes, as studied in \cite{ccfg,sym,racsam,cf-group}. While duality is a standard topic, inversion emerged through results from several papers:
\begin{itemize}
\item \cite{kaltpeck}, where it is constructed the twisted Hilbert space $Z_2$, who has a representation as $0\To \ell_2\To Z_2 \To \ell_2 \To 0$. In that same paper the authors show that $Z_2$ contains an Orlicz (non Hilbert) subspace $\ell_f$.
\item In \cite{cabe} it is shown that $Z_2/\ell_f$ is isomorphic to $\ell_f^*$, and therefore $Z_2$ has a representation of the form $0\To \ell_f\To Z_2 \To \ell_f^* \To 0$. The nature of $Z_2/\ell_f$ is not considered in \cite{kaltpeck}.
\item The general theory of twisted sums developed in \cite{kalt,kaltpeck} establishes a correspondence between exact sequences $0\To Y\To Z \To X \To 0$ of quasi Banach spaces and quasilinear maps $\Omega: X\To Y$. The quasilinear map generating
$0\To \ell_2\To Z_2 \To \ell_2 \To 0$ will be called $\KP$ in this paper. The quasilinear map generating $0\To \ell_f\To Z_2 \To \ell_f^* \To 0$ has been called $\mho$ in \cite{cabe}.
\item The connection between $\KP$ and $\mho$ was uncovered in  \cite{racsam} importing the notions of Domain and Range
from complex interpolation theory \cite{urbana}. Kalton  \cite{kaltdiff} had already established an intimate connection between complex interpolation and a special type of quasilinear maps called centralizers, paramount examples of which are $\KP$ and $\mho$. In \cite{racsam,sym} it was shown that if $\Omega_{\Psi, \Phi}$ is the differential generated by two interpolators $(\Psi, \Phi)$ then $ \Omega_{\Phi, \Psi}$ is a kind of inverse of $\Omega_{\Psi, \Phi}$. In particular, $\KP$ and $\mho$ will be inverses in the sense considered in this paper.
\end{itemize}

In this paper we transplant those ideas to a general context: each quasilinear map $\Omega: X\lop Y$ has an associated \emph{inverse} quasilinear map $\Omega^{-1}:\Ran(\Omega)\lop \Dom(\Omega)$, so that $\Dom(\Omega^{-1}) = Y$,  $\Ran(\Omega^{-1}) = X$ and both $\Omega^{-1}\Omega$ and $\Omega\Omega^{-1}$ are bounded (which somehow justifies the name).

Duality is another delicate point in the study of exact sequences of Banach spaces. It was considered from different angles in \cite{cabecastdu,castmoredu,cabeann,ccc}. In the second part of this paper we will first set the duality ideas in the broader context of this paper and study their connections with inversion. For instance, we will show that under reasonable conditions $(\Omega^{-1})^* = (\Omega^*)^{-1}$. The final Applications Section of the paper focuses on new results on duality and inversion for specific examples of quasilinear maps appearing in the literature.

\section{quasilinear maps and twisted sums in a broader context}\label{sect:ql-maps}

Let $X$ and $Y$ be quasi-Banach spaces with quasi-norms $\|\cdot\|_X$ and $\|\cdot\|_Y$.
We suppose that $Y$ is a subspace of some vector space $\Sigma$.

\begin{defin}\label{def:quasilinear}
A map $\Omega:X\to \Sigma$ is called \emph{quasilinear from $X$ to $Y$} with ambient space $\Sigma$ and denoted $\Omega: X\lop Y$ if it is homogeneous
and there exists a constant $C$ so that for each $x,z\in X$,$$
\Omega(x+z)-\Omega x-\Omega z\in Y\quad \textrm{and}\quad
\|\Omega(x+z)-\Omega x-\Omega z\|_Y \leq C(\|x\|_X +\|z\|_X).
$$
\end{defin}

The role of the ambient space was considered in \cite{cf-group}. There it was shown that given two quasilinear maps $\Phi, \Psi:X\lop Y$, the former with ambient space $A$ and the second with ambient space $B$ then one can consider that both maps have a certain common space $C$ as ambient space. The main implication of this is that $\Phi + \Psi$ exists.

\begin{defin}\label{def:bounded-q-l} A quasilinear map $\Omega:X\lop Y$ is said to be:
\begin{itemize}
\item \emph{bounded} if there exists a constant
$D$ so that $\Omega x \in Y$ and $\|\Omega x\|_Y\leq D\|x\|_X$ for each $x\in X$.
\item \emph{trivial} if there exists a linear map $L:X\To \Sigma$ so that $\Omega - L: X\To Y$ is bounded.
\end{itemize}\end{defin}

\begin{defin}\label{def:bounded-q-l} Let $\Phi, \Psi$ be quasilinear maps $X\lop Y$, that we can assume to have the same ambient space. They are said to be:
\begin{itemize}
\item \emph{boundedly equivalent} if $\Phi - \Psi$ is bounded.
\item \emph{equivalent} if $\Phi - \Psi$ is trivial.
\end{itemize}\end{defin}

A quasilinear map $\Omega: X\lop Y$ generates a quasi-Banach twisted sum of $Y$ and $X$ defined as
$$
Y\oplus_\Omega X =\{(\beta,x)\in \Sigma\times X : \beta-\Omega x\in Y\}
$$
and endowed with the quasi-norm $\|(\beta,x)\|_\Omega=\|\beta-\Omega x\|_Y+\|x\|_X$.
Indeed,
\begin{equation}\label{exact-seq1}
\begin{CD}0 @>>> Y @>\imath_1>> Y\oplus_\Omega X  @>\pi_2>> X@>>>0,\end{CD}
\end{equation}
where $\imath_1 (y) =(y,0)$ and $\pi_2(\beta,x)=x$ is an exact
sequence and $\|(y,0)\|_\Omega=\|y\|_Y$ and $\|(\Omega x,x)\|_\Omega=\|x\|_X$.

If $\Omega$ is bounded then $Y\oplus_\Omega X=Y\times X$ and $\|y-\Omega x\|_Y +\|x\|_X$
and $\|y\|_Y+\|x\|_X$ are equivalent quasi-norms on this space. If $\Omega$ is trivial then
$Y\oplus_\Omega X$ is isomorphic to $Y\oplus X$.\medskip

The map $B_\Omega: X\to Y\oplus_\Omega X$ defined by $B_\Omega x =(\Omega x,x)$
is a bounded homogeneous selection for $\pi_2$ such that
$\|B_\Omega x\|_\Omega= \|x\|_X$ for each $x\in X$. If $B_1: X\to Y\oplus_\Omega X$ is another homogeneous bounded selection for $\pi_2$
then it has the form $B_1 x = (\Omega_1 x,x)$ and $\Omega_1:X\lop Y$ is a quasilinear map boundedly equivalent to $\Omega$: indeed, given $x,z\in X$, since $B_1(x+z)-B_1x-B_1z =(\Omega_1 (x+z)- \Omega_1 x-\Omega_1 z,0)\in Y\oplus_\Omega X$,
we get $\Omega_1 (x+z)- \Omega_1 x-\Omega_1 z\in Y$ and
\begin{eqnarray*}
\|\Omega_1 (x+z)- \Omega_1 x-\Omega_1 z\|_Y &=&\|B_1(x+z)-B_1(x)-B_1z\|_\Omega\\
&\leq& C'\left(\|B_1(x+z)\|_\Omega +\|B_1(x+z)\|_\Omega\right) \leq C''(\|x\|_X +\|z\|_X).
\end{eqnarray*}
Moreover, for each $x\in X$, we have $(\Omega x,x), (\Omega_1 x,x)\in Y\oplus_\Omega X$,
hence $\Omega x-\Omega_1 x\in Y$ and\newline
$\|\Omega x-\Omega_1 x\|_Y= \|(\Omega x-\Omega_1 x,0)\|_\Omega\leq C(\Omega)(1+\|B_1\|)\|x\|_X$.\\

The notions of domain and range of centralizer-like mappings $\Omega: X\lop Y$ have been considered
in \cite{ccfg,racsam,sym,cf-group}.

\begin{defin}
Let $\Omega: X\lop Y$ be a quasilinear map with ambient space $\Sigma$.
\begin{itemize}
\item The domain of $\Omega$ is the space $\Dom\Omega=\{x\in X :\Omega x\in Y\}$ endowed with the quasi-norm $\|x\|_D =\|x\|_X+\|\Omega x\|_Y$.
Therefore $\Dom \Omega$ can be identified with the subspace $\{x\in X :(0,x)\in Y\oplus_\Omega X\}$ of $Y\oplus_\Omega X$.
\item The range of $\Omega$ is the space $\Ran\Omega=\{\beta\in \Sigma : \exists x\in X,\, \beta-\Omega x\in Y \}=
\{\beta\in\Sigma: \exists x\in X,\, (\beta,x)\in Y\oplus_\Omega X\}$, endowed with the quasi-norm
$\|\beta\|_R =\inf\{\|(\beta,x)\|_\Omega : (\beta,x)\in Y\oplus_\Omega X\}$.
Therefore $\Dom \Omega$ can be identified with the quotient space $\left(Y\oplus_\Omega X\right)/\Dom \Omega$.
\end{itemize}
\end{defin}

The domain and range of $\Omega$ are independent of the choice of the bounded selection chosen to define $\Omega$. More precisely:

\begin{lemma} Let $\Phi, \Psi$ be equivalent quasilinear maps $ X\lop Y$. Then
$\Dom \Phi = \Dom \Psi $ and $\Ran \Phi =\Ran \Psi $.\end{lemma}

Quite obviously,
\begin{prop}\label{prop:Om-bdd} $\Omega: X\To \Ran\Omega$ is bounded.\end{prop}
\begin{proof}
If $x\in X$ then $(\Omega x,x)\in Y\oplus_\Omega X$, hence $\Omega x\in\Ran\Omega$
and $\|\Omega x\|_R\leq \|(\Omega x,x)\|_\Omega =\|x\|_X$.
\end{proof}

If $\imath_2 x=(0,x)$ and $\pi_1(\beta,x)=\beta$ then one has an exact sequence
\begin{equation}\label{exact-seq2}
\begin{CD}0 @>>> \Dom\Omega @>\imath_2>> Y\oplus_\Omega X @>\pi_1>> \Ran\Omega@>>>0,\end{CD}
\end{equation}
Note that $\|(0,x)\|_\Omega=\|x\|_D$, $\|\beta\|_R\leq \|(\beta,x)\|_\Omega$ and
$\pi_1$ is surjective with $\ker \pi_1=\im \imath_2$.

\begin{corollary}\label{cor:dom-complete}
The spaces $(\Dom\Omega,\|\cdot\|_D)$ and $(\Ran\Omega,\|\cdot\|_R)$ are complete.
\end{corollary}
\begin{proof}
Since $\{x\in X :(0,x)\in Y\oplus_\Omega X\}$ coincides with $\ker \pi_1$, it is closed in
$Y\oplus_\Omega X$. A more pedestrian proof is possible: Suppose that $(0,x_n)$ converges to $(\beta,x)$ in $Y\oplus_\Omega X$.
Then $\beta\in \Ran\Omega$ and
$$
\lim_{n\to \infty}\|\beta-\Omega(x_n-x)\|_Y + \|x_n-x\|_X =0,
$$
hence $\beta-\Omega(x_n-x)\to 0$ in $Y$ and $x_n-x\to 0$ in $X$.
By Proposition \ref{prop:Om-bdd}, $\Omega(x_n-x)\to 0$ in $\Ran\Omega$.
Since $\|y\|_R\leq\|(y,0)\|_\Omega=\|y\|_Y$ for each $y\in Y$, we get $\beta=0$.\end{proof}

\section{Part I: Inversion}

We want to introduce the quasilinear map defining the sequence (\ref{exact-seq2}) and call it $\Omega^{-1}$. But, to preserve boundedness we need to be more specific:

\begin{prop}\label{prop:inv-Om-bdd} Let $\Omega:X\lop Y$ be a quasilinear map with ambient space $\Sigma$.
\begin{itemize}
\item  Let $M: \Ran \Omega \To X$ be an homogeneous map such that $B(\beta)= (\beta, M\beta)$ is a homogenous bounded selection $B: \Ran \Omega \To Y\oplus_\Omega X$ for $\pi_1$ in (\ref{exact-seq2}) with $\|B\|\leq K$.
\item Let $J: \Ran \Omega\to X$ be an homogeneous map such that $J(\omega) =x$ where $x\in X$ is such that $\|\omega  - \Omega x\| + \|x\| \leq K \mathrm{dist}((\omega,x), \Dom \Omega)$.
\end{itemize}
Then $M$ and $J$ are  boundedly equivalent quasilinear maps $\Ran \Omega \lop \Dom \Omega$.
\end{prop}
\begin{proof} It is clear that, by construction, both maps $M$ and $J$ define quasilinear maps $\Ran \Omega\lop \Dom \Omega $ with ambient space $X$. Indeed, given $\sigma, \mu \in \Ran \Omega$,
$M(\sigma+\mu)-M(\sigma) - M(\mu) \in \Dom \Omega$ if only we prove that $\Omega\left( M(\sigma+\mu)-M(\sigma) - M(\mu) \right)\in Y$. Since
$(\sigma, M\sigma), (\mu, M\mu), (\sigma + \mu, M(\sigma+ \mu)$ belong to $Y\oplus_\Omega X$ then also
$(0, M\sigma + M\mu - M(\sigma+ \mu)$ is in $Y\oplus_\Omega X$, which is what we want. Analogously with $J$.

To check their bounded equivalence, we show first that if $\omega\in \Ran \Omega$ then $M \omega - J\omega\in \Dom \Omega$. To get that it is enough to check that $(\omega, J\omega)\in Y\oplus_\Omega X$, as it clearly follows from $\|\omega  - \Omega J \omega \|\leq K\dist(\omega, \Dom \Omega)<\infty$, since then $(0, M\omega - J\omega)= (\omega, M\omega) - (\omega, J\omega) \in Y\oplus_\Omega X$. Finally

\begin{eqnarray*}
\|M \omega - J \omega\| &=& \|\Omega(M \omega - J \omega)\| \\
&=&\|\omega - \Omega(J\omega)  - \omega + \Omega(M\omega) - \Omega(M\omega) + \Omega(J\omega) + \Omega(M \omega - J \omega)\| \\
&\leq&\|\omega - \Omega(J\omega)\|  + \|\omega -\Omega(M\omega)\|  + \|M\omega\| + \|J\omega\| \\
&\leq& \|(\omega, M\omega)\|  + \|\omega -\Omega(J\omega)\|  + \|J\omega\| \\
&\leq& K\|\omega\| + \|(\omega, J\omega)|\|\\
&\leq& K\|\omega\| + K \mathrm{dist}((\omega, J\omega), \Dom \Omega)\\
&\leq& K\|\omega\| + K \|(\omega, J\omega) - (0, J\omega - M\omega)\\
&\leq& K\|\omega\| + K^2\|\omega\|.\qedhere
\end{eqnarray*}
\end{proof}

\adef We will call $\Omega^{-1}: \Ran \Omega \lop \Dom \Omega$ the quasilinear map in the previous proposition\zdef
In \cite{cabe}, the map $\Omega^{-1}$ was denoted $\mho$. Propositions \ref{prop:dom-ran} and \ref{prop:isomorphism} below show that $\Omega^{-1}$ can be seen as an inverse of $\Omega$. One has:
\begin{prop}\label{prop:dom-ran}$\;$
\begin{itemize}
\item $\Dom \Omega^{-1} =Y$
\item $\Ran \Omega^{-1} = X$
\item $\Omega^{-1}\Omega$ and $\Omega \Omega^{-1}$ are both bounded.
\item \cite[Proposition 3.8]{racsam}: In the particular case in which $\Omega$ is the differential $\Omega_{\Psi, \Phi}$ associated to a pair $(\Psi, \Phi)$ of interpolators on the same Kalton space then  $\Omega_{\Psi, \Phi}^{-1} = \Omega_{\Phi, \Psi}$, up to bounded equivalence.
\end{itemize}
\end{prop}
\begin{proof} If $\sigma \in \Ran \Omega \subset \Sigma$ is so that $\Omega^{-1} \sigma \in \Dom \Omega$ then there is $x\in X$ so that $\sigma - \Omega x\in Y$  and $\Omega x = \Omega \Omega^{-1} \sigma \in Y$. Consequently, $\sigma \in Y$. That $\Ran \Omega^{-1} = X$ is obvious. Indeed, an obvious choice for $\Omega^{-1}\Omega x$ is $x$. The other is analogous.\end{proof}
We can thus form the exact sequence
\begin{equation}\label{exact-seq3}
\begin{CD}0 @>>> \Dom\Omega @>\imath_1>> \Dom\Omega\oplus_{\Omega^{-1}} \Ran\Omega @>\pi_1>>
\Ran\Omega@>>>0,\end{CD}
\end{equation}
with $\imath_1x =(x,0)$ and $\pi_1(x,\beta)=\beta$ as usual; and show it is equivalent to (\ref{exact-seq2}):
The map $\Omega^{-1}$ defines the exact sequence (\ref{exact-seq2}) and, moreover:

\begin{prop}\label{prop:isomorphism}
There is a commutative diagram
$$\xymatrix{
&&Y\oplus_\Omega X\ar[dr]^{\pi_1} \ar[dd]^U&\\
0\ar[r]&\Dom\Omega \ar[ur]^{\imath_2}\ar[dr]^{\imath_1}& &  \Ran\Omega\ar[r]&0\\
&&\Dom\Omega \oplus_{\Omega^{-1}}\Ran\Omega \ar[ur]^{\pi_2}&}$$
\end{prop}
\begin{proof}
It is enough to observe that $U(\beta,x)=(x,\beta)$ is a bijective
operator that
makes commutative the diagram. Indeed, $U$ is well defined: if $(\beta,x)\in Y\oplus_\Omega X$ then $\beta\in\Ran\Omega$
and $(\beta,\Omega^{-1}\beta)\in Y\oplus_\Omega X$.
Thus $(0, x-\Omega^{-1}\beta)\in Y\oplus_\Omega X$, hence $x-\Omega^{-1}\beta\in\Dom\Omega$
and $(x,\beta)\in\Ran\Omega\oplus_{\Omega^{-1}} \Dom\Omega$.
\smallskip

$U$ is surjective: if $(x,z)\in\Ran\Omega\oplus_{\Omega^{-1}} \Dom\Omega$ then
$z\in\Ran\Omega$ and $x-\Omega^{-1}z\in \Dom\Omega$.
Thus $(z,\Omega^{-1}z), (0,x-\Omega^{-1}z)\in Y \oplus_\Omega X$, hence
$(z,x)\in Y\oplus_\Omega X$.
\smallskip

$U$ is bounded: $\|U(\beta,x)\|_{\Omega^{-1}}= \|(x,\beta)\|_{\Omega^{-1}}=
\|(x-\Omega^{-1}\beta\|_D +\|\beta\|_R$.
Since $\|\beta\|_R\leq \|(\beta,x)\|_\Omega$ and
\begin{eqnarray*}
\|x-\Omega^{-1}\beta\|_D &=&\|(0,x-\Omega^{-1}\beta)\|_\Omega =
\|(\beta,x)-(\beta,\Omega^{-1}\beta)\|_\Omega\\
&\leq& C\|(\beta,x)\|_\Omega+ C\|(\beta,\Omega^{-1}\beta)\|_\Omega
\leq C\|(\beta,x)\|_\Omega+ C\|B_p\| \|\beta\|_R,
\end{eqnarray*}
we get $\|U(\beta,x)\|_{\Omega^{-1}}\leq (C+C\|B_p\|+1)\|(\beta,x)\|_\Omega$, and
$U^{-1}$ is bounded by the open mapping theorem.
\end{proof}

With the same ideas we get that $\Omega$ and $(\Omega^{-1})^{-1}$ are boundedly equivalent. The dominion and range of a map can be used to detect its boundedness.

\begin{prop}\label{thm:bounded} A quasilinear map $\Omega: X \lop Y$ is bounded if and only if $\Dom\Omega=X$ and if and only if $\Ran\Omega=Y$. Therefore, $\Omega$ is bounded if and only if so is $\Omega^{-1}$. \end{prop}
\begin{proof}
The first assertion is immediate. Consequently, $\Omega^{-1}$ bounded if and only if $\Dom\Omega^{-1}=\Ran\Omega$.
Since $\Dom\Omega^{-1}=Y$ this is equivalent to $\Ran\Omega=Y$, which is equivalent to
$\Omega$ bounded. \end{proof}

Let us consider the diagram
\begin{equation}\label{sand-watch}
\xymatrixrowsep{1cm}
\xymatrixcolsep{2cm}
\xymatrix{& \Ran\Omega &\\
Y\ar[r]^{\imath_1} \ar[ru] & Y\oplus_\Omega X\ar[r]^{\pi_2}\ar[u]^{\pi_1} &X\\
 &\Dom\Omega \ar[u]^{\imath_2} \ar[ru]}
\end{equation}
the unlabeled arrows are the natural inclusions. The diagram is commutative since  $\pi_1 \imath_1 x =x $ and $\pi_2\imath_2 x=x$.

\begin{prop} Suppose that $\pi_2$ is strictly singular. Then $\pi_1$ is strictly singular if and only if the canonical inclusion $Y\To \Ran \Omega$ is strictly singular. Analogously, if $\pi_1$ is strictly singular then $\pi_2$ is strictly singular if and only if the canonical inclusion $\Dom \Omega \To X$ is strictly singular. \end{prop}
\begin{proof} Using the criterion \cite[Proposition 8]{cck} we get that since $\pi_2$ is strictly singular, $\pi_1$ is strictly singular if and only if so is $\pi_1\imath_1$.\end{proof}

\section{Part II. Duality}

To consider duality issues we must confine ourselves to deal with Banach spaces since quasi Banach spaces with trivial dual exist. But since a twisted sum of Banach space can be a quasi Banach (non Banach) space we must also restrict ourselves to consider only exact sequences $0\To Y \To Z \To X \To 0$ in which the three terms $Y,Z,X$ are Banach spaces. When $\Omega: X\lop Y$ is quasilinear then we only know that $Y\oplus_\Omega X$ is a quasi Banach space. To guarantee that it is isomorphic to a Banach space we must impose to the map $\Omega$ an additional condition:

\adef A quasilinear map $\Omega: X\lop Y$ is said to be a $1$-quasilinear map  if there is a constant $C$ such that for all finite sequences of elements $x_1, \dots, x_N \in X$
\begin{itemize}
\item[(a)] $\Omega(\sum_{n=1}^N x_n)- \sum_{n=1}^N \Omega(y_n)\in Y$
\item[(b)] $\|\Omega(\sum_{n=1}^N x_n)- \sum_{n=1}^N \Omega(x_n)\|_{Y}\leq C\sum_{n=1}^N \|x_n\|_X $ .
\end{itemize}
\zdef

If $\Omega$ is $1$-quasilinear then $\|\cdot\|_\Omega$ is equivalent to a norm, and thus $X\oplus_\Omega Y$ is a Banach space. Kalton showed \cite{kalt} that quasilinear maps on $B$-convex Banach spaces are $1$-quasilinear while Kalton and Roberts \cite{kaltrobe} showed that quasilinear maps on $\mathcal L_\infty$ spaces are $1$-quasilinear. We will say that a quasilinear map $\Omega:X\lop Y$ has \emph{dense domain}
if $\Dom\Omega$ is a dense subspace of $X$.

\begin{prop}\label{prop:dense} Let $\Omega:X \lop Y$ be a $1$-quasilinear map with dense domain. Then the map $J:\Dom\Omega\times Y\to Y\oplus_\Omega X$ given by $J(x,y)= (y,x)$ is
a continuous operator with dense range.
\end{prop}
\begin{proof}
Let $(x,y)\in\Dom\Omega\times Y$. Then
$$
\|J(x,y)\|_\Omega =\|y-\Omega x\|_Y+\|x\|_X\leq C(\|y\|_Y+\|\Omega x\|_Y)+\|x\|_X
\leq C(\|x\|_D+\|y\|_Y).
$$

For the second part, let $(\beta,x)\in Y\oplus_\Omega X$ and $\e>0$, and consider
the exact sequence (\ref{exact-seq1}).
Since $\Dom\Omega$ is dense in $X$, we can select $z\in\Dom\Omega$ such that
$\|q(\beta,x)-z\|_X<\e$.
Then $(\beta,x-z)\in Y\oplus_\Omega X$ and $\dist\left((\beta,x-z),\ker q\right)<\e$.
So we can choose $y\in Y$ such that $\|(\beta,x)-(y,z)\|<\e$.
\end{proof}

\begin{corollary}\label{cor:dense}
Let $\Omega:X \lop Y$ be a $1$-quasilinear map with dense domain. Then
\begin{enumerate}
  \item $\Omega^{-1}$ has dense domain; i.e., $Y$ is continuously embedded and dense in $\Ran\Omega$.
  \item $(\Ran\Omega)^*$ is continuously embedded in $Y^*$.
\end{enumerate}
\end{corollary}
\begin{proof} (1) If $y\in Y$ then $(y,0) \in Y\oplus_\Omega X$, hence $y\in \Ran\Omega$ and
$\|y\|_R\leq \|(y,0)\|_\Omega =\|y\|_Y$. Let $\beta\in\Ran\Omega$ and $\e>0$, and take $x\in X$ with $(\beta,x)\in Y\oplus_\Omega X$.
By Proposition \ref{prop:dense}, there is $(y,z)\in Y\times\Dom\Omega$ such that
$\|(\beta,x)-(y,z)\|_\Omega<\e$. Thus $\|\beta-y\|_R\leq \|(\beta,x)-(y,z)\|_\Omega<\e$. Note that $\Omega^{-1}$ is defined in $\Ran\Omega$ and, by Proposition \ref{prop:dom-ran}, $\Dom\Omega^{-1}=Y$. (2) follows from (1). \end{proof}

Consider the dual sequences of (\ref{exact-seq1}) and (\ref{exact-seq2}), namely:
\begin{equation}\label{exact-seq1d}
\begin{CD}0 @>>> X^* @>{\pi_2^*}>> (Y\oplus_\Omega X)^* @>{\imath_1^*}>> Y^*@>>>0,\end{CD}
\end{equation}
\begin{equation}\label{exact-seq2d}
\begin{CD}0 @>>> (\Ran\Omega)^* @>{\pi_1^*}>> (Y\oplus_\Omega X)^* @>{\imath_2^*}>>(\Dom\Omega)^*@>>>0.
\end{CD}
\end{equation}

By Proposition \ref{prop:dense}, the conjugate operator
$J^*:(Y\oplus_\Omega X)^*\to (\Dom\Omega)^*\times Y^*$ is continuous and injective.
Therefore, if $F\in (Y\oplus_\Omega X)^*$ then
$J^*F= (\alpha^*,y^*)\in (\Dom\Omega)^*\times Y^*$.

\begin{prop}\label{prop:J*}
Let $\Omega$ be a quasilinear map from $X$ to $Y$ with dense domain.
Let $F\in (Y\oplus_\Omega X)^*$ and suppose that $J^*F= (\alpha^*,y^*)$.
Then $j^* F =y^*$ and $i^*F=\alpha^*$.
\end{prop}
\begin{proof}
Let $F\in (Y\oplus_\Omega X)^*$, $y\in Y$ and $x\in\Dom\Omega$. Then

$\langle j^*F,y \rangle = \langle F,jy \rangle= \langle F,(y,0) \rangle=
\langle F,J(0,y)\rangle = \langle J^*F,(0,y)\rangle=
\langle(\alpha^*,y^*),(0,y)\rangle =\langle y^*,y\rangle$,
and similarly
$\langle i^*F,x \rangle =
\langle \alpha^*,x\rangle$.
\end{proof}

We pass to the definition of the dual $1$-quasilinear map.
Since $X^*$ is a subspace of $(\Dom\Omega)^*$, Proposition \ref{prop:J*}
allows us to define a map $\Omega^*$ as follows:

\begin{defin} Let $\Omega: X\lop Y$ be a $1$-quasilinear map with dense domain, and let $B$  be an homogeneous bounded selection for the quotient map $\imath_1^*$. We define the map $\Omega^*: Y^*\to (\Dom\Omega)^*$ by $J^*B y^*= (\Omega^*y^*,y^*)$. \end{defin}

\begin{prop}\label{thm:duality}
Let $\Omega: X\lop Y$ be a $1$-quasilinear map with dense domain. The map $\Omega^*: Y^* \lop X^*$  is $1$-quasilinear with ambient space  $(\Dom\Omega)^*$ and there exists an isomorphism from $X^*\oplus_{\Omega^*} Y^*$ onto $(Y\oplus_\Omega X)^*$.
\end{prop}
\begin{proof}
Since $B$ is a selector for $\imath_1^*$, for every finite set $y_1^*, \dots, y_n^*$ one has $\sum B(y_i^*)- B\left(\sum y_i^*\right) \in \ker \imath_1^*=\im \pi_2^*$ and thus
\begin{eqnarray*}
J^*\sum B(y_i^*)- J^*B\left(\sum y_i^*\right)&=& \sum \left( \Omega^*(y_i^*), y_i^*\right) - \left(\Omega^*\left(\sum y_i^*\right), \sum y_i^*\right)\\
&=&\left( \sum \Omega^*(y_i^*) - \Omega^*\left(\sum y_i^*\right),0 \right)\\
&\in& X^*\oplus_{\Omega^*} Y^*
\end{eqnarray*}
hence $\sum \Omega^*(y_i^*) - \Omega^*\left(\sum y_i^*\right)\in X^*$ and $
\|\sum \Omega^*(y_i^* )-\Omega^*\left(\sum y_i^*\right)\|_{X^*}\leq
2\|B\|\sum \|y_i^*\|_{Y^*}.$\\

We claim that $X^*\oplus_{\Omega^*} Y^*= J^*\left((Y\oplus_\Omega X)^*\right)$. Indeed, if $(\alpha^*,y^*)\in X^*\oplus_{\Omega^*} Y^*$ then $\alpha^*-\Omega^*y^*\in X^*$.
Thus $(\alpha^*-\Omega^*y^*,0)\in J^*\left((Y\oplus_\Omega X)^*\right)$, hence
$(\alpha^*,y^*) =(\alpha^*-\Omega^*y^*,0) +J^*By^* \in
J^*\left((Y\oplus_\Omega X)^*\right)$. Conversely, if $(\alpha^*,y^*)=J^* F$ for some $F\in (Y\oplus_\Omega X)^*$ then $\imath_1^*F=y^*$.
Thus $J^*(F-By^*)=(\alpha^*-\Omega^*y^*,0)\in \ker \imath_1^*$, hence
$\alpha^*-\Omega^*y^*\in X^*$, and we get $(\alpha^*,y^*)\in X^*\oplus_{\Omega^*} Y^*$.
\medskip

Now we consider the formal identity map $W:X^*\oplus_{\Omega^*} Y^*\to (\Dom\Omega)^*\times Y^*$ defined by
$W(\alpha^*,y^*)= (\alpha^*,y^*)$.
Since $W(\alpha^*,y^*)= (\alpha^*-\Omega^*y^*,0) +J^*B y^*$, we get
\begin{eqnarray*}
\|W(\alpha^*,y^*)\| &\leq& \|(\alpha^*-\Omega^*y^*,0)\| +\|J^*By^*\|\\
  &\leq& \|\alpha^*-\Omega^*y^*\|_{X^*}+\|J^*\|\cdot\|B\|\cdot \|y^*\|_{Y^*}\\
  &\leq& \|J^*\|\cdot\|B\|\cdot \|(\alpha^*,y^*)\|_{\Omega^*}.
\end{eqnarray*}
Thus $W$ is continuous.
Since $J^{*-1}$ is a closed operator, $V=J^{*-1}W$ is a closed bijective operator
from $X^*\oplus_{\Omega^*} Y^*$ onto $(Y\oplus_\Omega X)^*$.
By the closed graph theorem, $V$ is an isomorphism.
\end{proof}

We will refer to $\Omega^*$ as the \emph{adjoint of $\Omega$.} We have the following duality relations.
\begin{prop}\label{prop:dom-ran-duality}
Let $\Omega:X \lop Y$ be a $1$-quasilinear map with dense domain.
\begin{enumerate}
  \item $\Ran\Omega^* =(\Dom\Omega)^*$
  \item $\Dom\Omega^* =(\Ran\Omega)^*$.
\end{enumerate}
\end{prop}
\begin{proof}
(1) If $\alpha^*\in (\Dom\Omega)^*$, there is $F\in (Y\oplus_\Omega X)^*$ such that $\imath_1^*F=\alpha^*$.
Then $J^*F=(\alpha^*,y^*)$, hence $(\alpha^*,y^*)\in X^*\oplus_{\Omega^*} Y^*$, and
thus $\alpha^*\in \Ran\Omega^*$.

(2) Suppose that $y^*\in \Dom\Omega^* =\{y^*\in Y^* : \Omega^*y^*\in X^*\}$.
Then $(0,y^*)\in X^*\oplus_{\Omega^*} Y^*$, hence
$F=V(0,y^*)\in (Y\oplus_{\Omega} X)^*$. Since $i^* F=0$, we get $F\in \im p^*$,
hence $y^*\in (\Ran\Omega)^*$. If $y^*\in (\Ran\Omega)^*$, then $y^*\in Y^*$ by Corollary \ref{cor:dense}.
Moreover $V^{-1}p^*y^*=(0,y^*)\in X^*\oplus_{\Omega^*} Y^*$.
Then $\Omega^*y^*\in X^*$, hence $y^*\in \Dom\Omega^*$.
\end{proof}

The next result shows that $(\Omega^*)^{-1}$ and $(\Omega^{-1})^*$ are equivalent.

\begin{theorem} The exact sequence $$
\begin{CD}0 @>>> (\Ran\Omega)^* @>{\pi_1^*}>> (Y\oplus_\Omega X)^* @>{\imath_2^*}>> (\Dom\Omega)^*@>>>0,\end{CD}
$$ generated by the $1$-quasilinear map $(\Omega^{-1})^*$ and the exact sequence $$
\begin{CD}0 @>>> \Dom\Omega^* @>{\imath_2}>> X^*\oplus_{\Omega^*} Y^* @>{\pi_1}>> \Ran\Omega^*@>>>0\end{CD}$$
generated by $(\Omega^*)^{-1}$ are equivalent.
.\end{theorem}
\begin{proof} Proposition \ref{prop:dom-ran-duality} allows us to write the diagram
$$\begin{CD}
0 @>>>(\Ran\Omega)^* @>{\pi_1^*}>> (Y\oplus_\Omega X)^* @>{\imath_2^*}>> (\Dom\Omega)^* @>>>0\\
&& @| @V{J^*}VV @|\\
0 @>>>\Dom\Omega^* @>{\imath_2}>>  X^*\oplus_{\Omega^*} Y^* @>{\pi_1}>> \Ran\Omega^* @>>>0\\
\end{CD}$$
and Theorem \ref{thm:duality} shows that $J^*$ makes it commutative: $J^*\pi_1^* = (\pi_1 J)^*$ and $\pi_1 J^* = \imath_2^*$. \end{proof}

We describe next the duality action of $X^*\oplus_{\Omega^*} Y^*$ on (a dense subspace of) $Y\oplus_\Omega X$.

\begin{theorem}\label{prop:dual-action} Let $\Omega: X\lop Y$ be a $1$-quasilinear map with dense domain.
\begin{itemize}
\item Given $(\alpha^*,y^*)\in X^*\oplus_{\Omega^*} Y^*$ and $(y,z)\in J(\Dom\Omega\times Y)$,
\begin{equation*}
\langle V(\alpha^*,y^*),(y,z)\rangle
= \langle\alpha^* ,z \rangle + \langle y^*,y\rangle.
\end{equation*}
\item There exists $C>0$ such that for each $z\in Dom\Omega$ and $y^*\in Y^*$,
$$
\left|\langle y^*,\Omega z\rangle +\langle \Omega^*y^*,z\rangle \right|\leq C\|y^*\|_{Y^*}\|z\|_X.
$$
\end{itemize}
\end{theorem}
\begin{proof}
For the first part, $\langle V(\alpha^*,y^*),J(z,y)\rangle= 
\langle (\alpha^*,y^*),(z,y)\rangle = \langle\alpha^* ,z \rangle + \langle y^*,y\rangle$. For the second part, note that $(\Omega z,z)\in Y\oplus_\Omega X$, $(\Omega^* y^*,y^*)\in X^*\oplus_{\Omega^*} Y^*$
and
\begin{eqnarray*}
\left|\langle y^*,\Omega z\rangle+ \langle\Omega^* y^*,z\rangle \right|
&=& \left|\langle V(\Omega^* y^*,y^*),(\Omega z,z)\rangle\right|\\
&\leq& \|V\|\cdot\|(\Omega^* y^*,y^*)\|_{\Omega^*}\|(\Omega z,z)\|_\Omega\\
&\leq& \|V\|\cdot\|B_{j^*}\|\cdot\|y^*\|_{Y^*}\|z\|_X.\qedhere
\end{eqnarray*}
\end{proof}

Summing up since any exact sequence $0\To Y \To Z \To Z \To 0$ of Banach spaces has a dual sequence $0\To X^* \To Z^* \To X^* \To 0$,
one can associate a classical $1$-quasilinear map $\Omega^*: Y^*\To X^*$ to any classical $1$-quasilinear map $\Omega: X\To Y$. The quasilinear map $\Omega^*$ exists but is not unique. A different thing is what occurs with $1$-quasilinear maps $X\lop Y$ with ambient space $\Sigma$ as those we are dealing with in this paper. Let us first revisit the so-called Kalton duality explained in \cite{castmoredu} for classical quasilinear maps

\adef Two $1$-quasilinear maps $\Omega: B \to A$ and $\Phi: A^*\to B^*$ are called bounded duals one of the other if there is $C=C(\Omega, \Phi)>0$ such that for every $b\in B, a^*\in A^*$ one has
$$
|\langle \Omega b, a^*\rangle + \langle b, \Phi a^*\rangle |\leq C \|b\|\|a^*\|.$$\zdef

It is easy to check that if $\Phi, \Psi$ are both bounded duals of $\Omega$ then $\Phi - \Psi$ is bounded:
$$\left \langle (\Phi - \Psi)(x), y\rangle\right| = \left \langle \Phi x, y\rangle + \langle x, \Omega y \rangle - \langle x, \Omega y\rangle  - \langle \Psi x, y\rangle\right| \leq (C(\Phi, \Omega) +C(\Psi, \Omega)) \|x\|\|y\|$$

\begin{prop}\label{duality} If $\Omega$ generates the exact sequence $0\to A \stackrel{\imath}\to A\oplus_\Omega B \stackrel{\pi}\to B \to 0$ and $\Phi$ is a bounded dual of $\Omega$ then $\Phi$ generates the dual sequence $0\to B^*\to (A\oplus_{\Omega} B)^* \to A^* \to 0$, with the meaning that there is an operator
$D: B^*\oplus_\Phi A^* \longrightarrow (A\oplus_\Omega B)^*$ given by$$\langle D(b^*, a^*), (a,b)\rangle = \langle b^*, b\rangle + \langle a^*, a \rangle$$
making a commutative diagram
$$\xymatrix{
0\ar[r] &B^*\ar[r]^-{\imath}\ar@{=}[d] & B^*\oplus_\Phi A^* \ar[r]^-{\pi}\ar[d]^D & A^*\ar[r]\ar@{=}[d] &0\\
0\ar[r] &B^*\ar[r]_-{\pi^*} & (A\oplus_\Omega B)^* \ar[r]_-{\imath^*} & A^*\ar[r]&0}$$
In particular, the spaces $B^*\oplus_\Phi A^*$ and $(A\oplus_\Omega B)^*$ are isomorphic via the pairing
$$\langle (b^*, a^*), (a,b)\rangle = \langle b^*, b\rangle + \langle a^*, a \rangle$$
\end{prop}
\begin{proof} The continuity of $D$ is:
\begin{eqnarray*}
\langle b^*, b\rangle + \langle a^*, a \rangle &=&  \langle b^*, b\rangle + \langle a^*, a -\Omega b\rangle + \langle a^*, \Omega b\rangle \\
&=&  \langle b^* - \Phi a^*, b\rangle + \langle \Phi a^*,, b\rangle + \langle a^*, a -\Omega b\rangle + \langle a^*, \Omega b\rangle \\
&\leq& C|a^*\| \|b\| + \|a^*\|\|a- \Omega b\| + \|b^*\|\|b^* - \Phi a^*\|\\
&\leq& C'\|(a, b)\|_\Omega  \|(b^*, a^*)\|_\Phi.\end{eqnarray*}
Moreover $D$ makes the diagram commutative: $D\imath = \pi^*$ since
$$D\imath(b^*)(a,b) = D(b^*, 0)(a,b) = \langle b^*, b\rangle = b^*\pi(a,b) = \pi^*(b^*)(a,b)$$
while $\imath^* D = \pi$ since
$$\imath^* D(b^*, a^*)(a) = D(b^*, a^*)(a,0) = \langle a^*, a\rangle = \pi(b^*, a^*)(a).$$
Finally, $D$ is an isomorphism by the $3$-lemma.\end{proof}

The approach via Kalton duality presents some difficulties in the context of this paper in which quasilinear maps take values in some ambient space. Let us formulate the result in a typical situation: complex interpolation theory. Let $(X_0,X_1)$ be a Banach couple with sum (or ambient space) $\Sigma$ and intersection $\Delta$, which is equipped with the norm $x\in X_0\cap X_1\longmapsto \max\big{(} \|x\|_0, \|x\|_1\big{)}$. We assume that $(X_0,X_1)$ is  regular according to Cwikel \cite{cwdu}, i.e., $\Delta$ is dense in each $X_i$. Cwikel shows in \cite[Theorem 3.1]{cwdu}:\\

$\bigstar$ Each $X_i^*$ embeds into $\Delta^*$ in such a way that $X_0^*\cap X_1^*=\Sigma^*$.\\

Let $\Omega: X_0\cap X_1 \To \Sigma$ and $\Phi: X_0^*\cap X_1^* \To \Sigma$ be homogeneous maps. Observe that the duality conditions $\langle \Omega x, y\rangle$ and $\langle x, \Phi y\rangle$ make sense: $\Omega x\in \Sigma$, $y\in \Sigma^*$ and $x\in \Delta$, $\Phi y\in \Sigma =\Delta^*$.
Let now $X$ be a Banach space so that $\Delta \subset X\subset \Sigma$ and $\Sigma^*\subset X^* \subset \Delta^*$ with
$\Delta$ is still dense in $X$ and $\Sigma^* $ dense in $X^*$. Assume that $\Omega:X \lop X$ is a quasilinear map with ambient space $\Sigma$ and $\Phi: X^*\lop X^*$ is a quasilinear map with ambient space $\Delta^*$. If there exists $C>0$ such that for $x\in X_0\cap X_1=\Delta$, $y\in X_0^* \cap X_1^*=\Sigma^*$ one has $\left | \langle \Omega x, y\rangle + \langle x, \Phi y\rangle \right| \leq C \|x\|\|y\|$ then $\Omega, \Phi$ are bounded duals one of each other since there is an isomorphism $D(x^*,y^*)(y,x) = \langle x^*, x\rangle + \langle y^*, y\rangle$ making the diagram
$$\xymatrix{&&  X^*\oplus_\Phi X^*\ar[dd]^D\ar[dr] \\
0\ar[r]&X^* \ar[ur]\ar[dr] && X^* \ar[r]&0\\
&&(X\oplus_\Omega X)^* \ar[ur]& }$$
commute. This is a reformulation of what was already proved before and is precisely the context presented in \cite{ccc}. To provide a more general version we will use the following approach:

\begin{theorem}\label{dosestrellas} Let $X,Y$ be reflexive Banach spaces so that:
\begin{itemize}
\item Both $X$ and $X^*$ embed into an ambient space $\Sigma$ in such a way that there is a subspace $\Sigma_0$ that is dense both in $X$ and $X^*$ and $\Sigma\subset \Sigma_0^*$.
\item Both $Y$ and $Y^*$ embed into an ambient space $\Theta$ in such a way that there is a subspace $\Theta_0$ that is dense both in $Y$ and $Y^*$ and $\Theta\subset \Theta_0^*$.
\end{itemize}
Let $\Omega: X \lop Y$ be a $1$-quasilinear map with ambient space $ \Theta$ and such that $\Omega[\Sigma_0]\subset \Theta_0$ and let $\Phi: Y^*\lop X^*$ be a $1$-quasilinear map with ambient space $\Sigma$ and such that $\Phi[\Theta_0]\subset \Sigma_0$. If the following condition holds:

$\bigstar \bigstar$ There is $C>0$ so that for every $x\in \sigma$, $y\in \theta$ one has $\left | \langle \Omega x, y\rangle +  \langle x, \Phi y\rangle \right|\leq C \|x\|_X\|y\|_{Y^*}$

then $\Omega$ and $\Phi$ are bounded duals, so that the sequences $0\To Y \To Y\oplus_\Omega X \To X\To 0$ and $0\To X^* \To X^*\oplus_\Phi Y^* \To  Y^* \To 0$ are dual one of the other and the spaces $X^*\oplus_\Phi Y^*$ and $(Y\oplus_\Omega X)^*$ are isomorphic under the duality
$$ \langle (\eta, y^*), (\omega, x) \rangle = \langle \eta, x \rangle + \langle \omega, y^*\rangle.$$
\end{theorem}

We connect now inversion and duality.

\begin{lemma}\label{propperp} Under conditions Proposition \ref{dosestrellas}, $\left ( \Dom \Omega\right)^{\perp} = \Dom \Phi$ with the meaning $\mu \in \left ( \Dom \Omega\right)^{\perp} \Longleftrightarrow \jmath^* \mu \in \Dom \Phi$.\end{lemma}

\begin{proof} Pick $\mu\in (\Dom \Omega)^{\perp}$ with $\|\mu\|=1$. Since $\mu\in (\Dom \Omega)^{\perp} \Longleftrightarrow \mu(0,b)=0$ for all
$b\in \Dom \Omega$ then $\mu\left( (\Omega b, b)\right) = \mu\left( (\Omega b, 0)\right) = \langle \mu, \jmath\Omega b\rangle$ and also $|\langle \jmath^*\mu, \Omega b \rangle | = |\langle \mu, \jmath\Omega b\rangle\leq C \|b\|$ one gets from the estimate
$\left | \langle \Omega b, \jmath^*\mu\rangle + \langle b, \Phi \jmath^* \mu\rangle \right| \leq C \|b\|\|\mu\|$ that $\left | \langle b, \Phi \jmath^* \mu\rangle \right| \leq C' \|b\|$ for all $b\in \Dom \Omega$ and since $ \Dom \Omega$ is dense in $B$ by the hypotheses,
$\left | \langle b, \Phi \jmath^* \mu\rangle \right| \leq C' \|b\|$ for all $b\in B$; namely, $\Phi\jmath^* \mu \in B^*$ which means that $\jmath^* \mu\in \Dom \Phi$. \end{proof}

\begin{theorem} If $\Omega, \Phi$ are in the conditions of the previous proposition and satisfy condition $(\bigstar \bigstar)$ then $(\Phi)^{-1} \sim (\Omega^{-1})^*$.
\end{theorem}

\begin{proof} Since $\Omega$ and $\Phi$ are bounded duals, there is no loss of generality replacing $\Phi$ by $\Omega^*$. The proof is a drawing, in fact, two:
$$\xymatrixcolsep{2pc}
\xymatrix{&\Ran \Omega&\\
 \Dom \Omega^{-1}\ar[r]& Z \ar[r]\ar[u] & \Ran \Omega^{-1} \\
 &\Dom \Omega\ar[u]}\;\; \xymatrixcolsep{4pc}
\xymatrix{&(\Dom \Omega)^*&\\
 (\Ran \Omega^{-1})^*\ar[r]& Z^* \ar[r]\ar[u] & (\Dom \Omega^{-1})^* \\
 &(\Ran \Omega)^*\ar[u]}
$$
containing the two representations for $Z$ and their duals for $Z^*$. The exactness of the sequences for $Z$ means
$\Ran \Omega = Z/\Dom \Omega$, hence $(\Ran \Omega)^*=  (\Dom \Omega)^\perp$. Analogously
$(\Ran \Omega^{-1})^*=  (\Dom \Omega^{-1})^\perp$ and the right diagram becomes

$$\xymatrixcolsep{4pc}
\xymatrix{&(\Dom \Omega)^*&\\
 (\Dom \Omega^{-1})^\perp \ar[r]& Z^* \ar[r]\ar[u] & (\Dom \Omega^{-1})^* \\
 &(\Dom \Omega)^\perp\ar[u]}
$$
Since the horizontal sequence is the dual of the original sequence $\Omega$, and the vertical is the dual of $\Omega^{-1}$ the diagram is

$$\xymatrixcolsep{2pc}
\xymatrix{&\Ran ((\Omega^{-1})^{*})^{-1} = (\Dom \Omega)^*&\\
\Dom (\Omega^*)^{-1} = (\Dom \Omega^{-1})^\perp \ar[r]& Z^* \ar[r]\ar[u] & \Ran(\Omega^*)^{-1} = (\Dom \Omega^{-1})^* \\
 &\Dom ((\Omega^{-1})^{*})^{-1} =  (\Dom \Omega)^\perp\ar[u]}
$$

Apply now Proposition \ref{propperp}: use $(\Dom \Omega)^\perp = \Dom \Omega^*$ to get $\Dom (\Omega^*)^{-1} = (\Dom \Omega^{-1})^*$
and $(\Dom \Omega)^* = \Ran \Omega^*$ to get $\Ran(\Omega^*)^{-1}= (\Dom \Omega^{-1})^* = \Ran (\Omega^{-1})^*$. The horizontal sequence $\Omega^*$ is now
$$\xymatrix{
0\ar[r]&\Dom (\Omega^{-1})^{*}  \ar[r]& Z^* \ar[r] & \Ran (\Omega^{-1})^{*}\ar[r]&0}$$
therefore $\Omega^* = ((\Omega^{-1})^{*})^{-1}$ and thus $(\Omega^*)^{-1} = (\Omega^{-1})^{*}$.\end{proof}

\section{Applications}
Given an suitable pair of Banach spaces $(X_0, X_1)$ with ambient space $\Sigma$ one can consider the sequence of Schechter interpolators $\Delta_{k, \theta}(f) =f^{(k)}(\theta)(1/2)/k!$ defined on the suitable Calder\'on space $\mathcal C$ corresponding to the pair, see \cite{BL} or \cite{sche} for details. With those ingredients one can generate, following \cite{rochberg}, the family of associated Rochberg spaces
$\mathcal R_n(\theta) =\{ (\Delta_{n-1, \theta}(f), \dots, \Delta_{0, \theta}(f)): f\in \mathcal C\}$, that generalize the complex interpolation space $\mathcal R_1(\theta)=(X_0,X_1)_\theta$ and twisted sum space $\mathcal R_2(\theta)= \mathcal R_1(\theta)\oplus_{\Omega_\theta} \mathcal R_1(\theta)$, where $\Omega_\theta$ is the differential quasilinear map associated to the couple $(\Delta_1, \Delta_0)$. It was shown in \cite{cck} that for each pair $m,k\in\N$ with $n=m+k$ there are natural exact sequence
$$\begin{CD}
0@>>> \mathcal R_{m, \theta}  @>{\imath_{m,n}}>> \mathcal R_{n, \theta} @>{\pi_{n,k}}>> \mathcal R_{k, \theta} @>>> 0,
\end{CD}$$ that can be arranged forming commutative diagrams of exact sequences (we omit from now on the initial and final arrows $0\to \cdot$ and $\cdot \to 0$ and the parameter $\theta$)
\begin{equation*}\xymatrix{
&\mathcal R_k \ar@{=}[r]\ar[d]& \mathcal R_k \ar[d]&\\
&\mathcal R_n\ar[r]\ar[d]& \mathcal R_{n+m}\ar[r]\ar[d]& \mathcal R_m\\
&\mathcal R_{n-k}\ar[r]& \mathcal R_{n+m-k}\ar[r]& \mathcal R_m\ar@{=}[u]}
\end{equation*}
Let $\Omega_{m,n}: \mathcal R_m \lop \Sigma^n$ be the quasilinear map associated to the sequence $0\to \mathcal R_n\to \mathcal R_{n+m}\to \mathcal R_m \to 0$.

We shall consider three situations: first, the case of complex interpolation for the pair $(\ell_\infty, \ell_1)$ in which the involved differentials are not linear or bounded; second, translation operators which are bounded and homogeneous; third, the case of weighted K\"othe spaces, in which the differentials are linear unbounded operators.

\subsection{Higher order Kalton-Peck maps}\label{ex:1}  Our results will
complete those in \cite{sym}. The interpolators $\delta_{\theta}'$ and $\delta_{\theta}$ that define complex interpolation yield as interpolation space $X_\theta=(\ell_\infty, \ell_1)_\theta =\ell_p$ for $p=\theta^{-1}$ with associated differential $\KP: X_\theta \lop X_\theta$, usually called the Kalton-Peck map,
$$\KP (x)= p\; x\log \frac{|x|}{\|x\|_p}$$
and associated twisted sum space $X_\theta\oplus_{\KP} X_\theta = Z_p$, usually called the Kalton-Peck space \cite{kaltpeck}. According to \cite[Lemma 5.3 (c)]{kaltpeck}, see also \cite{sym}, $\Dom \KP$ is the Orlicz sequence space $\ell_{f_p}$ generated by $f_p(t)=t^p |\log t|^p$, while $\Ran \KP$ can be obtained by duality as the Orlicz space $\ell_{f_p}^* = \ell_{f_p^*}  $ generated by the Orlicz conjugate function $f_p^*$ of $f_p$ which, according to
\cite[Ex. 4.c.1]{lindtzaf}, is equivalent to $g_q(t)=t^q|\log t |^{p(1-q)} = t^q|\log t |^{-q}$ at $0$, for $pq=p+q$. We now follow \cite{cck,ccc,sym} to construct the associated Rochberg spaces. Let us focus at $\theta=1/2$: in this case each space $\mathcal R_n$ is  isomorphic to its dual and this allows us to compare the associated differential with its dual and uncover remarkable symmetries. Precisely, $\mathcal R_1=\ell_2$, while $\mathcal R_2 = Z_2$, the Kalton-Peck space. Let $\KP_{k,m}: \mathcal R_k\lop \mathcal R_m$ the $1$-quasilinear map with ambient space $\ell_\infty^m$ that generates the sequence $0\To \mathcal R_m  \To \mathcal R_n \To \mathcal R_k \To 0$. We are not interested here in the properties of the spaces $\mathcal R_n$, a topic studied in \cite{cck,ccc,sym}, but in the properties of the maps $\KP_{k,m}$. The case $m=k=1$ is in \cite{kaltpeck}. In general, one has:

\begin{theorem}
The $1$-quasilinear maps $\KP_{k,m}$ and ${(\KP_{m,k})}^*$ are isomorphic.
\end{theorem}
\begin{proof} It was proved in \cite{ccc} that for each $n\in\N$ there is an isomorphism
$u_n =\mathcal R_n\to \mathcal R_n^*$ given by
$$u_n\left(a_{n-i}\right)_{i=1}^{n} = \left((-1)^{n-i} a_{n-i}\right)_{i=1}^{n}.$$
such that for each $m,k\in\N$ and $n=m+k$ the following diagram is commutative:
\begin{equation}\label{eq:ZnZnd}
\begin{CD}
0@>>> \mathcal R_m  @>i_{m,n}>> \mathcal R_n @>\pi_{n,k}>> \mathcal R_k @>>> 0\\
&&@V(-1)^ku_mVV @V{u_n}VV @VV{u_k}V\\
0@>>> \mathcal R_m^* @>\pi_{n,m}^*>> \mathcal R_n^*  @>i_{k,n}^*>> \mathcal R_k^* @>>>0\\
\end{CD}\end{equation}
which implies that that ${\KP_{m,k}}^*$ is equivalent
to $(-1)^ku_m \KP_{k,m} u_k^{-1}$.\end{proof}

This yields the estimate:

\begin{eqnarray*}
|\langle \KP_{m,k}(a_0, \dots, a_{m-1}), (x_0, \dots, x_{k-1}) \rangle &+& \langle (-1)^ku_m \KP_{k,m} u_k^{-1}(x_0, \dots, x_{k-1}), (a_0, \dots, a_{m-1})\rangle|\\
&\leq& C \|(a_0, \dots, a_{m-1})\|_{\mathcal R_m} \|(x_0, \dots, x_{k-1})\|_{\mathcal R_k}\end{eqnarray*}

We know from \cite{cck} that \begin{eqnarray*}
\KP_{1,2}(x) &=& 2x(\log^2 \frac{|x|}{\|x\|}, \log \frac{|x|}{\|x\|})\\
\KP_{2,1}(y,x) &=& 2\left(y - 2x\log\frac{|x|}{\|x\|}\right) \log \frac{|y - 2x\log\frac{|x|}{\|x\|}|}{\|y - 2x\log\frac{|x|}{\|x\|}\|} + 2x\log^2\frac{|x|}{\|x\|}\end{eqnarray*}

and therefore the estimate $\left | \langle \KP_{1,2}(x), z\rangle - \langle x, u_1\KP_{2,1}u_2(z) \rangle\right| \leq C \|x\|_{\ell_2}\|z\|_{Z_2}$ becomes
\begin{eqnarray*}
&\;&\left |2xz_1\log^2 \frac{|x|}{\|x\|}+ 2xz_2 \log \frac{|x|}{\|x\|} - x \left(z_1 - 2z_2\log\frac{|z_2|}{\|z_2\|}\right) \log \frac{|z_1 + 2z_2\log\frac{|z_2|}{\|z_2\|}|}{\|z_1 + 2z_2\log\frac{|z_2|}{\|z_2\|}\|} - 2z_2\log^2\frac{|z_2|}{\|z_2\|} \right|\\ &\leq& C \|x\|_{\ell_2}\|(z_1, z_2)\|_{Z_2}
\end{eqnarray*}

Beware that the $u_n's$ cannot be deleted. For instance, $(\KP x, x)\in Z_2$ while $(\KP x, -x)\notin Z_2$, no matter if $Z_2$ and $u_2[Z_2]$ are isometric.

\subsection{Translation operators}\label{ex:ckmr}
Given a complex interpolation pair $(X_0, X_1)$ on the complex domain $U$ and $z, \theta \in U$, one can consider the pair of interpolators $(\Delta_{z}, \Delta_\theta)$ of evaluation at the points $z$ and $\theta$ respectively. Let $\texttt T_{z,\theta}$ be the associated quasilinear map, usually called the translation map. Let $B_\theta$ be a bounded homogeneous selector for $\Delta_\theta$, one has
$\texttt T_{z,\theta} = B_\theta(z)$ and therefore $\texttt T_{z,\theta}: X_\theta \To X_z$ is a bounded homogeneous map. Quite clearly one has
$\Dom \texttt T_{z, \theta}= X_\theta$ and $\Ran \texttt T_{z, \theta}= X_z$. Therefore $\texttt T_{z, \theta}^{-1}=\texttt T_{\theta, z}$, in accordance with the last entry in Proposition \ref{prop:dom-ran}. Translation operators appear also in the so-called differential methods of Cwikel, Kalton, Milman and Rochberg \cite{ckmr}, and the previous results remain valid in that context.

\subsection{Weighted K\"othe spaces}\label{ex:2}
Fix a K\"othe function space $X$ with the Radon-Nikodym property, let $w_0$ and $w_1$ be weight functions,
and consider the interpolation couple $(X_0, X_1)$, where $X_j= X(w_j)$, $j =0,1$ with their natural norms.
In \cite[Proposition 4.1]{ccfg} it was shown that $X_\theta = X(w_\theta)$ for $0<\theta <1$, where
$w_{\theta} = w_0^{1 - \theta}w_1^{\theta}$. In \cite{racsam} it was obtained
$X(\omega_0) \cap X(\omega_1) = X(\omega_\vee)$ with $\omega_\vee = \max\{\omega_0, \omega_1\}$;  $X(\omega_0) + X(\omega_1) = X(\omega_\wedge)$ with $\omega_\wedge = \min\{\omega_0, \omega_1\}$; $\Dom(\Omega_{\Psi, \Phi}) = X(w_\theta) \cap X(w_\theta \left|\log \frac{w_1}{w_0}\right|)$ and $\Ran(\Omega_{\Psi, \Phi})= X(w_\theta)+ X(w_\theta \left|\log \frac{w_1}{w_0}\right|^{-1})$. For $\Psi = \Delta'_{\theta}$ and $\Phi = \Delta_{\theta}$ we obtain
$\Omega_{\Psi, \Phi}f = \log \frac{w_1}{w_0}\cdot f$. Consequently
$\Omega_{\Psi, \Phi}^{-1}(x) = \Omega_{\Phi, \Psi}(x) = (\log \frac{\omega_1}{\omega_0})^{-1} x$. These maps are linear but unbounded, and in this case $\Dom \Omega_{\Psi, \Phi}f = \{f\in X(w): \log \frac{w_1}{w_0}\cdot f \in X(w)\} = X(w\log \frac{w_1}{w_0})$. Let us consider a specially interesting case already considered in \cite{sym}. Let $w$ be a weight sequence
(we will understand as in \cite[4.e.1]{lindtzaf} a non-increasing sequence of positive numbers such that $\lim w_n=0$ and $\sum w_n=\infty$). We set $w_0=w^{-1}$ and $w_1=w$ and let us consider the interpolation pair $(\ell_2(w^{-1}), \ell_2(w))$, whose complex interpolation space at $1/2$ is $\ell_2$. A homogeneous bounded selector for $\Delta_0$ is $B(x)(z) = w^{2z-1}x$. Thus $B(x)^{(n}(z) = 2^nw^{2z-1}\log^n w \cdot x$
and therefore $B(x)^{(n}(1/2) = 2^n\log^n w \cdot x$. Consequently,

$$\Omega_{\langle n-1, \dots, 1\rangle,0}(x) =  \left(\frac{2^{n-1}}{(n-1)!}\log^{n-1} w \cdot x, \dots, 2\log w \cdot x \right)$$

which yields
$$\mathcal R_n = \{(y_{n-1}, \dots, y_1, x)\in \ell_\infty^{n-1}\times \ell_2:  (y_{n-1}, \dots, y_1) - \Omega_{\langle n-1, \dots, 1\rangle,0}(x)\in \mathcal R_{n-1}\}$$
In particular, $\mathcal R_2 = \{(y,x)\in \ell_\infty\times \ell_2: x\in \ell_2, \quad y - 2\log w \cdot x\in \ell_2\}$ is generated by the differential $\Omega_{1,0}(x) = 2\log w \cdot x$. Consequently, $\Dom \Omega_{1,0} = \ell_2(\log w)$ and $\Ran \Omega_{1,0}=\ell_2(\log^{-1} w)$ and $\Omega_{1,0}^{-1}(x) = (2\log w)^{-1} \cdot x$. Analogously, it is now straightforward to check that
$$\Dom \Omega_{\langle n-1, \dots, 1\rangle,0} = \ell_2(\log^{n-1}w).$$


\begin{thebibliography}{99}



\bibitem{BL} J. Bergh and J. L\"ofstr\"om. \emph{Interpolation spaces. An introduction,}
Springer, 1976.

\bibitem{cabe} F. Cabello S\'anchez, \emph{Nonlinear centralizers with values in $L_0$,}
Nonlinear Anal. 88 (2014), 42--50.

\bibitem{cabeann} F. Cabello S\'anchez, \emph{Nonlinear centralizers in homology,} Math. Ann. 358 (2014), 779--798.

\bibitem{cabecastdu}
F.~Cabello S\'{a}nchez, J.M.F.~Castillo, \emph{Duality and
twisted sums of Banach spaces}, J. Funct. Anal. 175 (2000) 1--16.



\bibitem{ccc} F. Cabello S\'anchez, J.M.F. Castillo, W.H.G. Corr\^ea,
\emph{Higher order derivatives of analytic families of Banach spaces}, arXiv:1906.06677v2.

\bibitem{cck} F. Cabello S\'anchez, J.M.F. Castillo, N.J. Kalton,
\emph{Complex interpolation and twisted twisted Hilbert spaces},
Pacific J. Math. 276 (2015), 287--307.






\bibitem{ccfg} J.M.F. Castillo, W.H.G. Corr\^ea, V. Ferenczi, M. Gonz\'alez,
\emph{Stability properties of twisted sums generated by complex interpolation},
J. Inst. Math. Jussieu, to appear.  arXiv:1712.09647v3.



\bibitem{racsam} J.M.F. Castillo, W.H.G. Corr\^ea, V. Ferenczi, M. Gonz\'alez,
\emph{Differential processes generated by two interpolators}, RACSAM (2020) 114:183, 23 pp.  https://doi.org/10.1007/s13398-020-00920-5


\bibitem{sym} J.M.F. Castillo, W.H.G. Corr\^ea, V. Ferenczi, M. Gonz\'alez,
\emph{Interpolator symmetries and new Kalton-Peck spaces}, arXiv:2111.10640.

\bibitem{cf-group} J.M.F. Castillo, V. Ferenczi, \emph{Group actions on twisted sums of Banach spaces},  arXiv:2003.09767v2.




\bibitem{castmoredu} J.M.F. Castillo, Y. Moreno, \emph{Twisted dualities in Banach space theory},
in ``Banach spaces and their applications in analysis", in Honor of Nigel Kalton's 60th Birthday,
pp. 59--76. (B. Randrianantoanina and N. Randrianantoanina eds.), Walter de Gruyter, 2007,






\bibitem{cwdu} M. Cwikel, \emph{Lecture Notes on Duality and Interpolation Spaces}, (2014), arXiv:0803.3558v2, 36 pp.


\bibitem{urbana} M. Cwikel, B. Jawerth, M. Milman, R. Rochberg,
\emph{Differential estimates and commutators in interpolation theory,}
in ``Analysis at Urbana II'', London Mathematical Society, Lecture Notes Series 138,
E.R. Berkson, N.T. Peck, J.J. Uhl eds., pp. 170--220, Cambridge Univ. Press, 1989.

\bibitem{ckmr}  M. Cwikel, N.J. Kalton, M. Milman, R. Rochberg,
\emph{A unified theory of Commutator Estimates for a class of interpolation methods},
Adv. in Math. 169 (2002) 241-312.








\bibitem{kalt} N.J. Kalton, \emph{The three-space problem for locally bounded F-spaces}, Compositio Math. 37 (1978) 243--276.

\bibitem{kaltdiff} N.J. Kalton,
\emph{Differentials of complex interpolation processes for K\"othe function spaces,}
Trans. Amer. Math. Soc. 333 (1992), 479--529.


\bibitem{kaltpeck} N.J. Kalton, N.T. Peck,
\emph{Twisted sums of sequence spaces and the three space problem,}
Trans. Amer. Math. Soc. 255 (1979), 1--30.


\bibitem{kaltrobe} N.J. Kalton, J.W. Roberts, \emph{Uniformly exhaustive submeasures and
nearly additive set functions},  Trans. Amer. Math. Soc. 278 (1983) 803--816.


\bibitem{lindtzaf} J. Lindenstrauss, L. Tzafriri, \emph{Classical Banach spaces I, sequence spaces}, Springer, 1977.





\bibitem{rochberg} R. Rochberg,
\emph{Higher order estimates in complex interpolation theory,}
Pacific J. Math. 174 (1996), 247--267.




\bibitem{sche} M. Schechter, \emph{Complex interpolation,} Compositio Math. 18 (1967), 117--147.


\end{thebibliography}
\end{document}